\documentclass{amsart}

\usepackage{amssymb}
\usepackage[alphabetic,initials]{amsrefs}

\usepackage{graphicx}

\newtheorem{theorem}{Theorem}[section]

\newtheorem{corollary}[theorem]{Corollary}

\theoremstyle{definition}

\theoremstyle{remark}
\newtheorem{remark}[theorem]{Remark}

\numberwithin{equation}{section}

\usepackage{setspace}

\begin{document}

\title[The Random Graph and the Curve Graph]{The Random Graph Embeds in the Curve Graph of any Infinite Genus Surface}
\author{Edgar A. Bering IV}
\address{Department of Mathematics, Statistics, and Computer Science \\
University of Illinois at Chicago \\
851 S Morgan St. (M/C 249) \\
Chicago, IL \\
60607 }
\email{eberin2@uic.edu}

\author{Jonah Gaster}

\address{Department of Mathematics, Boston College \\ Chestnut Hill, MA 02467}
\email{gaster@bc.edu}


\date{}

\begin{abstract}
The random graph is an infinite graph with the universal property that any
embedding of $G-v$ extends to an embedding of $G$, for any finite graph.  
In this paper we show that this graph embeds in the curve graph of 
a surface $\Sigma$ if and only if $\Sigma$ has infinite genus, showing that the
curve system on an infinite genus surface is ``as complicated as possible''. 
\end{abstract}

\maketitle

\section{Introduction}

In this paper we will prove

\begin{theorem}\label{thm:maintheorem}
	The random graph embeds into the curve graph $\mathcal{C}(\Sigma)$ if
	and only if $\Sigma$ has infinite genus.
\end{theorem}

We adopt the terminology that an \emph{embedding} of a graph $f: G\to H$ is a
one-to-one map on the vertices so that $(u,v)$ is an edge in $G$ if and only if
$(f(u),f(v))$ is an edge in $H$. (This is also called an \emph{induced
subgraph} elsewhere in the literature.)

The one-ended, infinite genus, orientable surface with one boundary component
is a subsurface of any orientable infinite genus surface~\cite{infinitegenus}.
The choice of a disk on $\Sigma$, the one-ended orientable surface of infinite
genus without boundary, thus induces an embedding of the curve graph
$\mathcal{C}(\Sigma)$ into the curve graph of an arbitrary orientable infinite
genus surface. 
Therefore, for one direction of the theorem, it suffices to produce an
embedding of the random graph into $\mathcal{C}(\Sigma)$ when $\Sigma$ is the
one-ended orientable surface of infinite genus.

The other direction is 
perhaps more surprising, since it is tempting to view an infinite-type surface of finite genus as already quite complicated. 
However, Ehrlich, Even, and Tarjan~\cite{eet} showed that there are graphs that cannot be realized as the incidence graph of a collection of planar intervals (in their language, there are graphs not of \emph{planar type}), and we employ
their construction to demonstrate the necessity of infinite genus. 

Rado~\cite{rado} showed that every countable graph embeds in the random graph,
however we focus on the random graph for its combinatorial properties. The
first order theory of the random graph, in the graph language, is not edge
stable in the model-theoretic sense~\cite{unstable}. It follows that

\begin{corollary}
	The first order theory of the curve graph of an infinite genus surface
	is not edge stable.
\end{corollary}

The lack of edge stability implies that the theory of the curve graph is also
unstable in the model theoretic sense. With Gabriel Conant, we prove a
complementary result for finite-type surfaces~\cite{bcg}; the theory of
$\mathcal{C}(\Sigma)$ for a finite-type surface $\Sigma$ is edge stable. The
combination of these two results show that the model-theoretic dividing line of
edge stability and the topological dividing line of finite-type coincide for
curve graphs. It is unknown whether or not the curve graph of a finite-type
surface is stable.

Erd\H{o}s and R{\'e}nyi introduced the random graph from a probabilistic point
of view, constructing a graph on countably many vertices by adding edges with
probability $\frac{1}{2}$. The result of this construction is almost surely
isomorphic to a unique object, which we call \emph{the random
graph}~\cite{erdosrenyi}. Rado gave an explicit construction of the random
graph: take as vertices the natural numbers $\mathbb{N}$. Given $x < y$, add an
edge $(x,y)$ if the $x$th bit of the binary expansion of $y$ is 1~\cite{rado}.
The random graph enjoys a universal property, known as the extension property;
for any finite graph $G$, if $G-v$ embeds into the random graph then this
embedding can be extended to $G$. 

The other graph of interest in this article is the curve graph of an infinite
genus surface, with or without boundary. A simple closed curve on a surface is
essential if no component of the complement is a disk, and non-peripheral if no
component of the complement is an annulus. For brevity, we will use curve to
mean the isotopy class of an essential non-peripheral simple closed curve. The
intersection number of two curves (denoted $i(\alpha,\beta)$) is the minimum
cardinality of the intersection taken over all transverse realizations of $\alpha$ and
$\beta$.
 
Fix a surface $\Sigma$. The curve graph $\mathcal{C}(\Sigma)$ has as vertices
the curves on $\Sigma$, and an edge between the vertices corresponding to
curves $\alpha$ and $\beta$ when $i(\alpha,\beta) = 0$.  (As an aside, this
construction can be extended to a definition of a higher-dimensional simplicial
complex of interest in its own right, but we will focus on the
1-skeleton~\cites{harvey,masurminsky1,masurminsky2}.) When $\Sigma$ is of
finite-type, $\mathcal{C}(\Sigma)$ is well-known to be $\delta$-hyperbolic and
infinite-diameter, with automorphism group isomorphic to the mapping class
group of $\Sigma$~\cites{ivanov,masurminsky1}.  Recent work has focused on
different analogues of $\mathcal{C}(\Sigma)$ when $\Sigma$ is of infinite-type,
more satisfying from a geometric viewpoint~\cites{bavard,
aramayona-fossas-parlier, aramayona-valdez, fossas-parlier,
durham-fanoni-vlamis}.

In fact, it is not hard to see that every finite graph $G$ embeds in the curve graph of a surface for some
closed surface of genus $g$. We outline a simple proof. Suppose $G$ has $n$
vertices.  Let $\Sigma_{0}$ indicate a closed surface of genus large enough so
that $\Sigma_{0}$ contains a collection of $n$ curves in minimal position that
pairwise intersect once\footnote{Genus $\lceil \frac{n-1}{2} \rceil$
suffices. Such a system of curves has been referred to as a \emph{complete
1-system} in the literature.}, and identify these curves with the vertices of
$G$ arbitrarily. For each edge between a pair of vertices of $G$, add a handle
near the intersection of the corresponding curves on $\Sigma_{0}$, and thread
one of the curves through the handle so that the new curves on the new surface
do not intersect. The identification of the vertices of $G$ with the resulting
curve system extends to an embedding of $G$ into the curve graph
$\mathcal{C}(\Sigma)$ of the resulting surface $\Sigma$, by construction. 

\begin{remark}
This leaves open the problem of determining the minimal genus such that
every finite graph on $n$ vertices embeds in the curve graph of that
genus
(cf.~\cite{kim-koberda}*{Question 1.1}).
The above construction provides the bound $O(n^2)$.
\end{remark}

This implies that every finite graph embeds into the curve graph of an infinite
genus surface, and it suggests that this curve graph of an infinite genus
surface should be quite complicated. Note, however, that this alone does not
guarantee the presence of the random graph. Also note that it is apparent that
the random graph does not embed in the curve graph of any finite type surface
surface: A complete subgraph of the curve graph of a surface of genus $g$ with
$n$ punctures has size at most $3g-3+n$, whereas every finite complete graph
embeds in the random graph. 

\begin{figure}
\begin{center}
\includegraphics{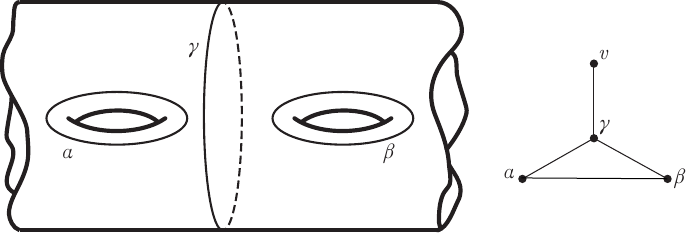}
\caption{A collection of curves and a partial graph embedding that cannot extend.}\label{fig:nonisom}
\end{center}
\end{figure}

Moreover, the curve graph of an infinite genus surface cannot itself
be isomorphic to the random graph. Fix an infinite genus surface $\Sigma$. Pick
a separating curve\footnote{Such a curve always exists: If $\gamma$ is
non-separating, choose a curve that intersects it once, and take a regular
neighborhood of the union. The boundary of this neighborhood is a separating
curve.} 
$\gamma$ and two curves $\alpha,\beta$, one in each component of
$\Sigma \setminus \gamma$. Let $G$ be the graph in figure \ref{fig:nonisom}. We
can embed $G-v$ into $\mathcal{C}(\Sigma)$ according to the labeling in the
figure, but an extension to $v$ would imply the existence of a curve disjoint
from $\gamma$ that intersects both $\alpha$ and $\beta$, a contradiction since
$\gamma$ is separating.

\section{Proof of Theorem \ref{thm:maintheorem}}

\begin{proof}

We deal with the forward implication first, showing that the presence of the
random graph in the curve graph implies that $\Sigma$ has infinite
genus.  It is evidently enough to show that there is a finite graph
which does not embed in $\mathcal{C}(\Sigma)$ when $\Sigma$ has finite
genus, since any finite graph embeds in the random graph.  This
construction is due to Ehrlich, Even, and Tarjan when $g=0$~\cite{eet}; 
it is simple enough to include completely. 

Suppose that $\Sigma$ has genus $g< \infty$, and choose a finite graph $G$ that
admits no topological embedding $G \to \Sigma_{g,0}$ into the closed surface of genus $g$ (a large complete graph
will do). Consider the graph $G'$ obtained by adding a vertex that
subdivides each edge of $G$, so that there are now $|V(G)|$ \emph{old}
vertices and $|E(G)|$ \emph{new} vertices of $G'$.  Let $\hat{G}$ indicate
the complementary graph of $G'$, and note that: (1) each new vertex $v$ of
$\hat{G}$ is incident to all other vertices of $\hat{G}$, except for the two old
vertices that are incident to the edge of $G$ corresponding to $v$, and (2) each old vertex is incident to every other old vertex.

Suppose that $\hat{G}$ is realized by a curve system $\Gamma$ on $\Sigma$ in
minimal position. The subdivision of the vertices of $\hat{G}$ into new
and old vertices gives a subdivision of $\Gamma$ into new and old
curves.  For each old curve $\gamma$, select a point
$p_{\gamma}\in\gamma$ in the complement of its intersections with the
other curves of $\Gamma$, and contract $\gamma\setminus p_{\gamma}$ to
a point.  Because the old vertices of $\hat{G}$ are all incident to
each other, when we do this contraction to each old curve one-by-one,
we obtain $|V(G)|$ on vertices on $\Sigma_{g,0}$.  Moreover, because
the new vertices of $\hat{G}$ are all incident to each other, the new
curves become a system of disjoint arcs connecting these vertices.  By
construction the resulting arcs provide a topological embedding of $G$
into $\Sigma_{g,0}$, a contradiction.

For the other direction, 
we will provide an explicit construction of a family of curves on the one-ended
orientable infinite genus surface whose intersection relation is exactly
described by the random graph. Our approach is in two parts; first we will give
a countable collection of multicurves with this property, then describe how to
add handles to convert these multicurves to curves without changing the
intersection relation of the curve system or the homeomorphism type of the surface.

Rado's construction fits more naturally into the setting of multicurves, so we
first define a multicurve complex $m\mathcal{C}(\Sigma)$ analogous to the curve
complex. Let the vertices of $m\mathcal{C}(\Sigma)$ be finite sets of disjoint
curves. For multicurves $U,V\in m\mathcal{C}(\Sigma)$, let $i(U,V)$ be the sum of
intersection numbers $i(\alpha,\beta)$ over all $\alpha\in U,\beta\in V$. In
analogy with the curve graph, there is an edge in $m\mathcal{C}(\Sigma)$ between
$U,V$ if $i(U,V)=0$. 
(The multicurve graph has also been called the \emph{clique graph} in the literature~\cite{kim-koberda}).
Below, we write multicurves additively,
e.g. $\alpha+\beta$ is the multicurve $\{\alpha,\beta\}$.

To fix notation, let $\Sigma_1$ be the one-ended orientable surface of infinite
genus. Note that the random graph is self-complementary (that is, the
complement graph is isomorphic to the random graph), so we will work with the
complement of Rado's model: let $x,y\in\mathbb{N}$ with $x<y$ be adjacent when
the $x$th bit in the binary expansion of $y$ is $0$. We describe below a map
$[\cdot]:\mathbb{N}\to m\mathcal{C}(\Sigma_1)$ so that, for $x < y$, the
intersection number $i([x],[y])$ is equal to the $x$th bit in the binary
expansion of $y$. Such a map induces an embedding of the random graph into
$m\mathcal{C}(\Sigma_1)$. 

Realize $\Sigma_1$ in $\mathbb{R}^3$ as the regular neighborhood of the lattice
on points $\mathbb{N}\times\{0,1\}\times\{0\}$. With this embedding the
`centers' of `holes' of $\Sigma_1$ occur at $(n+\frac{1}{2},\frac{1}{2},0)$ with
$n\in\mathbb{N}$. The intersection of $\Sigma_1$ with the coordinate plane
$\mathbb{R}\times\mathbb{R}\times\{0\}$ is the disjoint union of countably many
circles and a real line. Let $a_{i}$ be the circle component in the strip
$(i-1,i)\times\mathbb{R}\times\{0\}$, and let $b_i$ be the Dehn twist of $a_i$
around the intersection of the half-plane
$\{i-\frac{1}{2}\}\times(-\infty,\frac{1}{2})\times\mathbb{R}$ with $\Sigma$.
In other words, $a_{i}$ winds around the $i$th hole of $\Sigma$, and
$i(a_i,b_j) = \delta_{i,j}$, as pictured in Figure \ref{fig:aandbcurves}. 

\begin{figure}
\begin{center}
\includegraphics{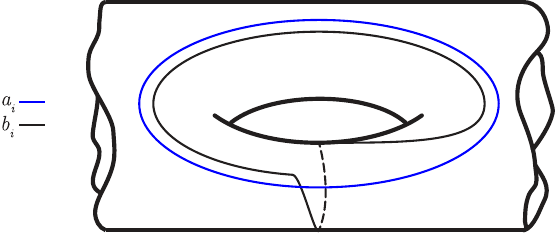}
\caption{A possible choice for $a_i$ and $b_i$ at handle $i$.}\label{fig:aandbcurves}
\end{center}
\end{figure}

Given a natural number $x$, let $x_i$ be the $i$th binary digit in the
expansion of $x$. We define
\begin{equation*}
[x] = b_x+\sum_{i=0}^{\lceil\log_2 x\rceil} x_i\cdot a_i.
\end{equation*}
Figure $\ref{fig:homologyexample}$ shows $[0]$ and $[4]$.  By construction this
is our desired map and the intersection relation among the multicurves
$\{[n]\}_{n\in\mathbb{N}}$ is given by the random graph.

\begin{figure}
\begin{center}
\includegraphics{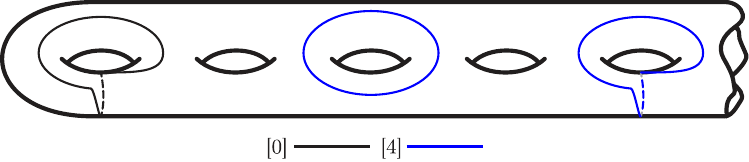}
\caption{An illustration of geometric realizations of $[0]$ and $[4]$.}\label{fig:homologyexample}
\end{center}
\end{figure}

\begin{figure}
\begin{center}
\includegraphics{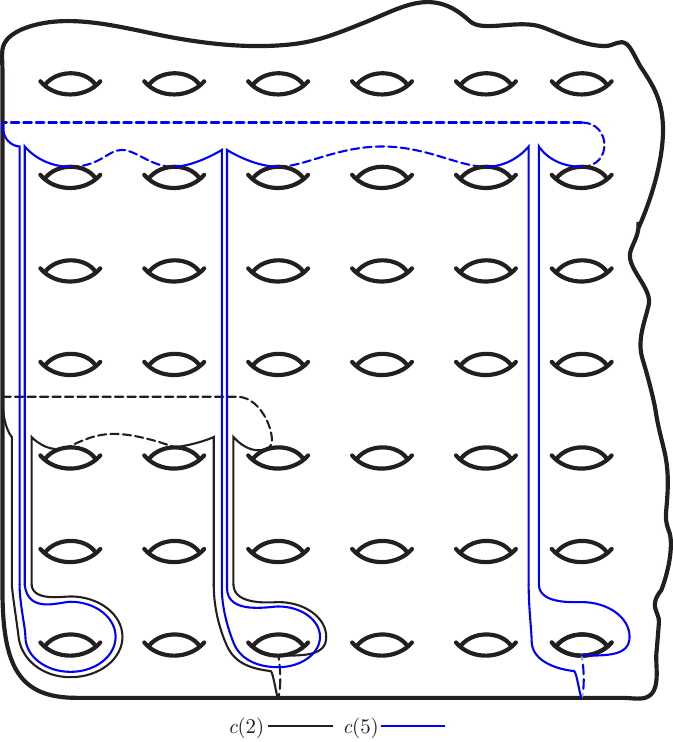}
\caption{An illustration of $c(2)$ and $c(5)$ realizing $i(c(2),c(5)) = 1$.}\label{fig:curveex}
\end{center}
\end{figure}

At this point one would like to blindly add handles to realize these
multicurves as curves. However, for each bit there are infinitely many curves
that need to be connected to the handle encoding that bit, so care must be
exercised. Consider a new realization of $\Sigma_1$ in $\mathbb{R}^3$, as the
regular neighborhood of the lattice on $\mathbb{N}\times\mathbb{N}\times\{0\}$.
The centers of `holes' are now at $(x+\frac{1}{2},y+\frac{1}{2},0)$ for $x,y\in
\mathbb{N}^2_{\geq 0}$. This naturally indexes the rows and columns of the
embedding (row $n$ is the regular neighborhood of points of the form $(x,n,0)$,
and the columns are similarly indexed).  We take $a_i,b_j$ as before (along the
$x$-axis). For a multicurve $[x] = b_x+\sum_{i} x_i\cdot a_i$, construct the
curve $c(x)$ by connecting each $a_i$ or $b_i$ in $[x]$ to row $x+1$ by a pair
of vertical lines along column $i$, and then connect these arcs to one another
along the `back' of $\Sigma_1$; figure \ref{fig:curveex} shows $c(2)$ and $c(5)$.
For $x < y$, we can realize $c(x)$ and $c(y)$ so that when $x$ and $y$ use a
common column $c(x)$ passes outside of $c(y)$; hence $i(c(x),c(y)) =
i([x],[y])$. (Note that, when curves intersect once, this intersection is
necessarily essential \cite{farb-margalit}.) We conclude that
$\{c(n)\}_{n\in\mathbb{N}}$ is the vertex set of an embedding of the random
graph in $\mathcal{C}(\Sigma_1)$.

\end{proof}

\begin{remark}
	The embedding of the random graph given by $\{c(n)\}_{n\in\mathbb{N}}$
	is far from unique. Each curve $b_i$ is constructed by a single Dehn
	twist. Varying the powers of each twist defining a $b_i$ individually
	yields systems of curves in distinct mapping class group orbits. The
	extended mapping class group of $\Sigma$ is isomorphic to the graph
	automorphisms of $\mathcal{C}(\Sigma)$ in the case of the infinite
	genus surface with one end, so these embeddings are also
	combinatorially inequivalent~\cite{hernandez-valdez}.
\end{remark}

\section*{Acknowledgments}

The first author thanks Marc Culler for helpful conversations. Both authors
thank Tarik Aougab for communicating interesting questions about the complexity
of curve graphs, and Gabriel Conant for pointing out the model-theoretic
corollary.

\end{document}